\input ./style/arxiv-general.cfg
\documentclass[aop,MSNbibl,seceqn,dvips]{arximspdf}
\makeatletter
   \@ifpackageloaded{graphicx}{}{\usepackage{graphicx}}
\makeatother

\doi{10.1214/15-AOP1001}
\volume{44}
\issue{2}
\pubyear{2016}
\firstpage{1285}
\lastpage{1307}
\docsubty{FLA}

\makeatletter
\def\sfrac#1#2{#1/#2}
\def\vfrac#1#2{(#1)/#2}

\def\sklfrac#1#2{(#1/#2)}

\newcommand{\rrvert}{\vert}
\newcommand{\rrVert}{\Vert}
\newcommand{\llvert}{\vert}
\newcommand{\llVert}{\Vert}
\renewcommand{\mid}{|}
\newtheorem{theorem}{Theorem}[section]
\newtheorem{lemma}[theorem]{Lemma}
\newproclaim{remark}[theorem]{Remark}
\newtheorem{proposition}[theorem]{Proposition}
\newproclaim{question}[theorem]{Question}
\newproclaim{definition}{Definition}
\newcommand{\R}{\mathbb{R}}
\newcommand{\Z}{\mathbb{Z}}
\newcommand{\N}{\mathbb{N}}
\def\dd{d}
\def\diam{\operatorname{diam}}
\def\P{\mathbf{P}}
\def\E{\mathbf{E}}
\def\diam{\operatorname{diam}}
\def\supp{\operatorname{supp}}
\def\auxi{\hat{\operatorname{aux}}}
\def\Aux{\operatorname{Aux}}
\def\D{{\mathcal D}}
\def\exp{\operatorname{exp}}
\makeatother

\begin{document}
\begin{frontmatter}

\title{A Poisson allocation of optimal tail}
\runtitle{A Poisson allocation of optimal tail}

\begin{aug}
\author[A]{\fnms{Roland}~\snm{Mark\'o}\ead[label=e1]{roland.marko@hcm.uni-bonn.de}\thanksref{T1}}
\and
\author[B]{\fnms{\'Ad\'am}~\snm{Tim\'ar}\corref{}\ead[label=e2]{madaramit@gmail.com}\thanksref{T2}}
\runauthor{R. Mark\'o and A. Tim\'ar}
\affiliation{Universit\"{a}t Bonn and Alfr\'ed R\'enyi Institute of Mathematics}
\address[A]{Hausdorff Center for Mathematics\\
Universit\"{a}t Bonn\\
D-53115 Bonn\\
Germany\\
\printead{e1}}
\address[B]{Alfr\'ed R\'enyi Institute of Mathematics\\
Re\'altanuda u. 13-15\\
H-1053 Budapest\\
Hungary\\
\printead{e2}}
%
\end{aug}
\thankstext{T1}{Supported in part by a Hausdorff scholarship.}
\thankstext{T2}{Supported by MTA R\'enyi ``Lend\"ulet'' Groups and
Graphs Research Group and by TAMOP 4.2.4.
A/1-11-1-2012-0001 National Excellence Program, a project subsidized by
the European
Union and co-financed by the European Social Fund. The author is a Marie Curie Fellow.}

%
\received{\smonth{3} \syear{2013}}
%
\revised{\smonth{12} \syear{2014}}

%
\begin{abstract}
The allocation problem for a $d$-dimensional Poisson point process is
to find a way to
partition the space to parts of equal size, and to assign the parts to the
configuration points in a measurable, ``deterministic'' (equivariant) way.
The goal is to make the diameter $R$ of the part assigned to a configuration
point have fast decay. We present an algorithm for $d\geq3$ that
achieves an $O(\exp(-cR^d))$ tail, which is optimal up to $c$. This improves
the best previously known allocation rule, the gravitational allocation,
which has an $\exp(-R ^{1+o(1)})$ tail. The construction is based on the
Ajtai--Koml\'os--Tusn\'ady algorithm and uses the
Gale--Shapley--Hoffman--Holroyd--Peres stable marriage scheme (as
applied to
allocation problems).
\end{abstract}

%
\begin{keyword}[class=AMS]
\kwd{60D05}
\end{keyword}
\begin{keyword}
\kwd{Fair allocation}
\kwd{Poisson process}
\kwd{translation-equivariant mapping}
\end{keyword}
\end{frontmatter}

\section{Introduction}\label{s.intro}

Consider the random discrete point set $\omega$ in $\R^d$ given by the
Poisson point process of intensity 1. We would like to find functions
$\psi_\omega\dvtx  \omega\to L^1(\R^d)$ that assign to each point of
$\omega
$ a set of measure 1, and such that $\psi_\omega$ is a measurable,
equivariant function of $\omega$. Then we call $\omega\mapsto\psi
_\omega$ an allocation rule. See Definition~\ref{def2} for more details.

We prove the existence of an allocation rule of the following (optimal)
tail.

\begin{theorem} \label{main}
For $d\geq3$ there exist $c,b>0$ and an allocation rule $\omega
\mapsto\psi_\omega$ for the
Poisson point process such that
\[
\P\bigl[\diam\bigl(\{0\}\cup\psi_\omega(0)\bigr)>R \mid 0\in \omega
\bigr] \leq c\exp \bigl(-bR^d\bigr).
\]
\end{theorem}

The proof of Theorem~\ref{main} is constructive. The allocation rule
presented here is built upon a generalization of the algorithm due to
Ajtai, Koml\'os and Tusn\'ady \cite{AKT} and employs a local variant of
the stable marriage allocation introduced by Hoffman, Holroyd and Peres
\cite{HHP}, based on the ``stable marriage'' algorithm of Gale and Shapley.

Informally, in this allocation problem, we want to divide land between
a set of farmers randomly scattered in space, and in such a way that
each farmer
knows the rule to determine his own land up to a small error, by
looking at the locations of other farmers in a big
enough neighborhood. This rule has to be the same for everybody.
We want to find an allocation rule where the maximal distance a farmer
has to walk (fly) between any two points of his land is minimal. That is,
we want the tail
\[
\P\bigl[\diam\bigl(\{0\}\cup\psi_\omega(0)\bigr)>R \mid 0\in\omega\bigr]
\]
to decay as fast in $R$ as possible, where $\psi_\omega(0)$ denotes
the cell of $0$
(conditioned on having a center in $0$).

Now we give the precise definitions.
In our previous, less formal definition, we wanted to partition the
space and assign pieces to the centers. Here it will be more convenient
to use the indicator functions of these pieces.
Let $\Omega$ be the set of discrete sets of
points in $ \R^d$. The points of an $\omega\in\Omega$ are called
\textit{centers}. Let us fix an $\omega\in\Omega$.

\begin{definition}[(Allocation)]\label{def1}
An \textit{allocation} is a function $\psi_\omega\dvtx \omega\to L^1(\R^d)$
such that:
\begin{longlist}[(2)]
\item[(1)]
for every $\xi\in\omega$, $\psi_\omega(\xi)$ is a function with
values from $\{0,1\}$;
\item[(2)] for Lebesgue almost every $x\in\R^d$, there is at most one
$\xi\in\omega$ such that $\psi_\omega(\xi)(x)=1$;
\item[(3)]
for every $\xi\in\omega$, $\int_{\R^d} \psi_\omega(\xi)(x) \,dx=1$.
\end{longlist}

Without assumption $(3)$ we call $\psi_\omega$ a \textit{weak allocation}.
We say
that $\psi_\omega\dvtx \omega\to L^1(\R^d)$ is a \textit{weak fractional
allocation} if
$(1')$ and $(2')$, below, hold and a \textit{fractional allocation} if
$(1')$, $(2')$ and $(3)$ hold:
\begin{longlist}[(2$'$)]
\item[(1$'$)]
for every $\xi\in\omega$, $\psi_\omega(\xi)$ is a function with
values from $[0,1]$;
\item[(2$'$)]
for Lebesgue almost every $x\in\R^d$, $\sum_{\xi\in\omega} \psi
_\omega
(\xi)(x)\leq1$.
\end{longlist}
\end{definition}

Sometimes it will be natural to think about a (fractional) allocation
$\psi_\omega\dvtx  \omega\to L^1 (\R^d)$ as a family $\{\psi_\omega
(\xi)
\dvtx \xi\in\omega\}$ of
functions.

\begin{definition}[(Allocation scheme, allocation rule)]\label{def2}
Let $\P$ be the law of the Poisson point process of intensity 1 in $\R
^d$. An
\textit{allocation scheme} is a mapping $\omega\mapsto\psi_\omega$
that is
defined for $\P$-almost every $\omega\in\Omega$, measurable (i.e., if
$(S, \mathcal S, \mu)$ is the underlying probability space of the
Poisson point process, and $L^1(\R^d)$ is equipped with the Borel sets
induced by the $L^1$-distance, the mapping $s \mapsto(\omega(s),\psi
_{\omega(s)})$ from $S \to\{(\omega, \phi) \mid  \omega\in\Omega,
\phi
\in[L^1 (\R^d)]^\omega\}$ is measurable with respect to $\mathcal
S$), and such that almost surely
$\psi_\omega$ is an allocation. If,\vspace*{1pt} furthermore, the mapping $\omega
\mapsto\psi_\omega$ is translation-equivariant; that is, for any
$\omega\in\Omega$, $\xi\in\omega$ and $x,y\in\R^d$ we
have
\[
\psi_{\omega+x}(\xi+x) (y+x)=\psi_{\omega}(\xi) (y),
\]
then we call the allocation scheme an \textit{allocation rule}.
Define \textit{weak}, \textit{fractional} and \textit{weak fractional allocation
rules} (\textit{schemes}) analogously.
\end{definition}

Some use the term \textit{allocation factors} of the point process
for allocation rules. We mention that allocation rules satisfy
something stronger than (2) from Definition~\ref{def1}; namely, for
Lebesgue almost every $x\in\R^d$, there is \textit{exactly one}
$\xi\in\omega$ such that $\psi_\omega(\xi)(x)=1$; see, for example,
\cite{HP} for a proof of this simple statement.

Denote by $0$ the origin in $\R^d$ and by $\lambda$ the Lebesgue
measure in $\R^d$. The letters $b$ and $c$ will stand for positive real
constants, whose value may change when employed in different
statements. From now on, \textit{we always assume that} $d\geq3$.

We want to define an allocation rule in such a way that the probability
$\P[\diam(\psi_{\omega\cup\{0\}} (0)\cup\{0\} )\geq R  \mid    0\in
\omega]$ decays as fast as possible.
By translation invariance we
could have taken any other fixed $\xi_0 \in\R^d$ instead of $0$, and
the tail
would be the same.

Define $\Omega':=\{\omega\cup\{0\}\dvtx   \omega\in\Omega\}$ and
$\P'
=[\cdot   \mid 0\in\omega]$
as the Palm version of the probability measure $\P$ that defines the
Poisson point process (and $\E'$ the
corresponding expectation). To facilitate readability we will tend to
use the notation $\omega'$ for elements of $\Omega'$. It is well known
that $\P'
[\{\omega\cup\{0\}\dvtx   \omega\in E\}]=\P[E]$ for every measurable $E$
in $\Omega$; see, for example, \cite{SKM} for the proof of this
statement and for other basic facts and definitions about the Palm
version of point processes. Rephrasing the previous paragraph, our goal
is to make $\P'[\diam(\psi_{\omega'}(0) \cup\{0\})>R]$ have fast decay.

In several papers (\cite{HHP} and follow-up works) a slightly different
problem is the focus of interest. Suppose we have the allocation rules
$\omega\mapsto\psi_\omega$ and $\omega\mapsto\phi_\omega$. The
rule $\omega\mapsto\psi_\omega$ defines a unique center $\xi_0=\xi
_0(\omega, \psi_\omega)$ with $0 \in\psi_\omega(\xi_0)$ almost
surely. Then
%
\begin{equation}
\label{coupling} \P\bigl[\diam\bigl(\phi_{\omega-\xi_0}(0)\cup\{0\}\bigr)>R\bigr]=
\P'\bigl[\diam \bigl(\phi_{\omega'}(0)\cup\{0\}\bigr)>R\bigr],
\end{equation}
where the equation follows from Theorem 13 in \cite{HP} (the claim that
$\xi_0$ is a so-called nonrandomized extra head scheme).

The objective in the setup of the cited papers is again to obtain a
rule with optimal tail bound, this time for the random variable $\llvert  \xi
_0\rrvert  $. Our setup is stronger than this one, meaning that for any
allocation rule $\omega\mapsto\psi_\omega$ the tail probability $\P
[\diam(\{0\}\cup\psi_\omega(0))>R   \mid  0\in\omega]$ is greater
than or equal to $\P[|\xi_0(\omega, \psi_\omega)|  > R]$, hence any
upper bound on the decay of the former (as in Theorem~\ref{main})
implies an upper bound on the decay of the latter (as in \cite{HHP}).
To see this, set $\phi=\psi$ and apply~(\ref{coupling}), so that we have
%
\begin{eqnarray}
\label{coupling2} \P\bigl[\llvert \xi_0\rrvert > R\bigr] &\leq& \P\Bigl[
\max_{x \in\psi_\omega(\xi_0)} \llvert x-\xi_0\rrvert >R\Bigr]\leq\P
\bigl[\diam\bigl(\psi_\omega(\xi_0)\cup\{\xi_0\}
\bigr)>R\bigr]
\nonumber
\\
&=&\P\bigl[\diam\bigl(\psi_{\omega-\xi_0}(0)\cup\{0\}\bigr)>R\bigr]
\\
&=& \P'\bigl[\diam\bigl(\psi_{\omega'}(0)\cup\{0\}\bigr)>R\bigr].\nonumber
\end{eqnarray}
In fact, the current setup is strictly stronger: imagine an allocation
rule when most cells consist of ball-like pieces of almost unit volume
containing the corresponding center, and a small extra piece far away
from this one. For such an allocation, the quantity on the left of
(\ref{coupling}) would decay quickly, while on the right we could have slow decay.

The relation between the two discussed tail events allows us to phrase
a lower bound on $\P'[\diam(\psi_{\omega'}(0)\cup\{0\})>R]$ which
every allocation rule has to satisfy. Clearly,
%
\begin{equation}
\qquad \P\bigl[\llvert \xi_0\rrvert > R\bigr] \geq\P\bigl[B(0,R) \cap
\omega= \varnothing\bigr]=\exp \bigl(-\lambda \bigl(B(0,R)\bigr)\bigr)=\exp
\bigl(-cR^d\bigr),
\end{equation}
where $B(0,R)$ is the Euclidean ball around $0$ with radius $R$, and
$c>0$ is a constant only depending on the dimension $d$. Therefore it
follows through (\ref{coupling2}) for any allocation rule $\omega
\mapsto\psi_\omega$ that
\[
\exp\bigl(-cR^d\bigr) \leq\P'\bigl[\diam\bigl(
\psi_{\omega'}(0)\cup\{0\}\bigr)>R\bigr],
\]
which explains the claim regarding the optimality of our construction.

The allocation problem was first studied in a finite setup, where finitely
many points are distributed uniformly and independently in a box.
Here, of course, the requirement of equivariance is meaningless. We will
present later a variant of the algorithm by Ajtai, Koml\'os and Tusn\'ady \cite{AKT},
which was a crucial component of several later methods for the finite problem.
For $n$ uniformly independently distributed points in a cube of volume
$n$, it was proved \cite{AKT} that the average diameter of
an
allocation cell is $\log^{1/2} n$ for $d=2$ and finite for $d\geq3$,
and precise rates of decay were determined subsequently; see \cite{TY} for
details and the sharpest results.
Interest in the
infinite setup originated from the fact that an allocation rule gives
rise to a
shift-coupling between a
point process and its Palm version; see, for example, \cite{HL,HP}. In~\cite{HP}, Holroyd and Peres studied the problem of how to
find the optimal tail of an allocation rule. In the same paper, they
presented a \textit{randomized} invariant allocation rule of optimal tail
decay. (Randomized allocation rules are defined similarly to allocation
rules in Definition~\ref{def2}, but the use of extra randomness is
allowed; i.e., the allocation scheme is not necessarily a deterministic
function of the point configuration.) Let us mention that several
related optimization problems are much easier to handle for the
randomized variant, such as in the case of matching schemes of point
processes or coin flips. Deterministic constructions have been an area
of active research ever since, and allocation rules that satisfy
additional conditions have been subjects of recent analysis (e.g.,
stability \cite{HHP}, connectedness \cite{K}). The
best previously known allocation rule for $d\geq3$ was the \textit{gravitational
allocation}, investigated in \cite{CPPR,CPPR2}, where the tail
is $\P
[\llvert  \xi_0(\omega)\rrvert  >R]= \exp(-R^{1+o (1)})$. In \cite{HP}, the assumption
$d\geq3$ is necessary for an
exponential tail, where $\E[\llvert  \xi_0(\omega)\rrvert  ^{d/2}]=\infty$ is
proved for
$d=1,2$. The best currently known upper bounds for the tails in the
case of $d=1$ and $d=2$ were\vspace*{1pt} presented in \cite{HPPS}, showing $\P
[\llvert  \xi
_0(\omega)\rrvert   > R]\leq c R^{-1/2}$ and $\P[\llvert  \xi_0(\omega)\rrvert   > R]\leq c
R^{-0.496\ldots}$, respectively, and they were achieved by the stable
allocation of Hoffman, Holroyd and Peres \cite{HHP}. Briefly, the reason
for the drastic change of behavior from dimension 3 is
that here the isoperimetric function of $\R^d$ becomes larger in magnitude
than the deviation of the number of points in a ball.
For the existence of optimal allocations
with respect to other quantities (e.g., the \textit{average} distance of a
center from the points of the cell) and connections to optimal transport,
see \cite{HS}.

\subsection{Construction}\label{s.cons}
\begin{figure}

\includegraphics{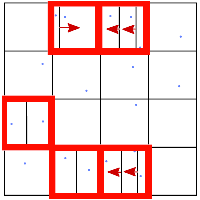}

\footnotesize{(a) The initial stage with centers}\vspace*{9pt}

\includegraphics{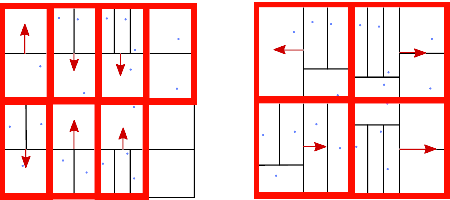}

\footnotesize{(b) The two steps of the first stage}\vspace*{9pt}

\includegraphics{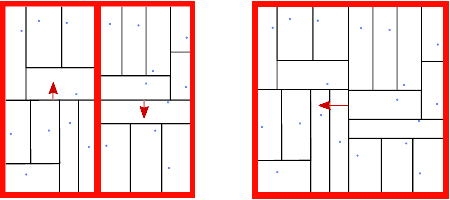}

\footnotesize{(c) The two steps of the second stage}\vspace*{9pt}

\includegraphics{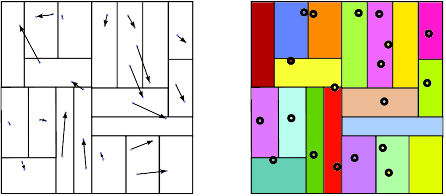}

\footnotesize{(d) The position of centers and final auxiliary points, the resulting local allocation}

\caption{AKT local allocation between $v+[0,2^2)^d$ and the centers in
$d=2$ ($v$ is the bottom left corner). Transformations take place
inside the highlighted cuboids; the bisector walls touched by the
arrows are moved in the indicated direction. (The length of an arrow
may not be proportional to the length of the shift, for better perspicuity.)}\label{abra}
\end{figure}

Let us present the (surprisingly simple) construction for Theorem~\ref
{main} briefly, before going
into the details. First, for any $v\in\R^d$ we will define a sequence
of weak allocation schemes
associated to it. These will not be equivariant yet. The construction
will be based on a straightforward generalization of the algorithm of
Ajtai, Koml\'os and Tusn\'ady (\textit{AKT algorithm}) \cite{AKT}, which
assigns a piece of unit volume to each of $n$ points in a box $B$, in
such a way
that the pieces partition $B$.
Furthermore, if the points are scattered
uniformly and independently (which is the same as the restriction of
the Poisson point process to $B$,
conditioned on there being $n$ points in $B$), then the average
diameter of the cells has asymptotically the same tail behavior as in
Theorem~\ref{main} (if the volume of $B$ and the number of random
points in it are asymptotically the same).
We extend the method to the Poisson point process in $\R^d$ by
subdividing $\R^d$ to the cubes of size $2^n$ in $v+2^n\Z^d+[0,2^n)^d$,
and applying the AKT algorithm to each of the cubes. The result is a
weak allocation scheme for each $n$ and $v$, which we call $\mathrm{AKT}_{\omega,v,n}$, and denote the cell of $\xi\in\omega$ under this
scheme by $f_{v,n}^{\omega,\xi}$.
The algorithm for given $n$ will be called the AKT($v$) algorithm run
up to stage $n$. The details of the AKT algorithm are discussed in
Section~\ref{s.akt}. See also Figure~\ref{abra}.

Having defined for every $v\in\R^d$ a sequence of weak allocation
schemes dependent on $v$, we next want to remove this dependence and
construct a sequence of weak allocations whose elements are equivariant
(i.e., weak allocation rules).
We will see that for every $\xi\in\omega$
%
\begin{equation}
f_{v,n}^{\omega,\xi}=f_{u,n}^{\omega,\xi}\qquad\mbox{whenever }
u-v\in2^n\Z^d. \label{e-1}
\end{equation}
Now, define
%
\begin{equation}
f_n^{\omega,\xi}:=\frac{1}{2^{nd}} \int_{[0,2^n)^d}
f_{v,n} ^{\omega
,\xi}\,\dd v, \label{f_n}
\end{equation}
a function from $\R^d$ to $[0,1]$. This is well defined (the integral
exists) by Lemma~\ref{Lemma2} below.
It is equivariant by (\ref{e-1}), so it is a fractional weak allocation
rule for each $n \in\N$. Now we want to get rid of ``weakness.''
We will prove that with probability~1 all the $f_n^{\omega,\xi}$ have an
$L^1$ limit ($\xi\in\omega$),
%
\begin{equation}
f^{\omega,\xi}:=\lim_{n \to\infty} f_n^{\omega,\xi},
\label{f_infty}
\end{equation}
and that it is a function of integral 1. From the
properties of the AKT construction (Lemma~\ref{fontos} below) we will show
that the diameter of the support of this limit function has the tail
that we want. We will
conclude that the map $\eta\dvtx  \omega\mapsto f^\omega$ with $f^\omega
\dvtx \xi\mapsto f^{\omega,\xi}$ ($\xi\in
\omega\in\Omega$) defines a fractional allocation rule with the desired
tail.

Finally, we will define an allocation $\omega\mapsto\psi_\omega$ with
$\psi_\omega\dvtx \xi\mapsto\psi_\omega(\xi)$ ($\xi\in
\omega$) from the above fractional
allocation. It will be such that the support of $\psi_\omega(\xi)$ is
contained in the support
of $f^{\omega,\xi}$. This step (Lemma~\ref{final}) will be
based on the fact that almost every point of $\R^d$ is contained in the
support of $f^{\omega, \xi}$
for only finitely many $\xi$.

In the next section we give a summary of the AKT construction, presenting
the generalized version that we are using. Section~\ref{rule} continues
with the sequence of invariant weak allocations, the
limiting fractional allocation and finally, the allocation rule that
satisfies Theorem~\ref{main}. This section ends with some concluding
remarks. In Section~\ref{technical} we give the necessary bounds for
the concentration of the cell diameter, which come by technical
modifications of the similar, usual bounds for the AKT method.

\section{The generalized AKT algorithm, bounds}\label{s.akt}

The AKT method was developed by Ajtai, Koml\'os and Tusn\'ady in \cite
{AKT}, and it outputs a local allocation (see definition below) between
a finite cube and the i.i.d. uniform random points lying in it. We
require the following notion.

\begin{definition}
Let $C$ be a measurable bounded subset of $\R^d$ and $\omega_C$ a
finite, nonempty subset of $C$. Call $\psi\dvtx  \omega_C \to L^1(C)$ a
\emph{local allocation} between $C$ and $\omega_C$ if:
\begin{longlist}[(2)]
\item[(1)]
for every $\xi\in\omega_C$, $\psi(\xi)$ takes values from $\{0,1\}$;
\item[(2)] for Lebesgue almost every $x\in C$ there is exactly one
$\xi\in\omega_C$ such that $\psi(\xi)(x)=1$;

\item[(3)]
for every $\xi\in\omega_C$, $\int_{C} \psi(\xi)(x) \,\dd x=\frac
{1}{\llvert  \omega_C\rrvert  }\lambda(C)$.
\end{longlist}
\end{definition}

%
Fix $\omega\in\Omega$, and assume for simplicity that for any two
points of $\omega$, their respective coordinates are pairwise different
(which has probability 1 if they are distributed according to the
Poisson point process).




We now present the local allocation given by the AKT scheme between the
cube $C=v+[0,2^N)^d$ and $\omega_C:=\omega\cap C$, where $v \in\R^d$
and $N \geq1$. See Figure~\ref{abra} for an illustration.
We define the method recursively with respect to $N$.

For $N=0$ and $\llvert  \omega_C\rrvert   =k >0$, let us fix some arbitrary way of
dividing a unit cube containing $k$ points into $k$ connected parts
with equal measure and assigning them to the points $\xi$. (For $k=0$,
we do not do anything.)
For example,
separate the centers by $k-1$ hyperplanes orthogonal to the first axis
such that the $i$th hyperplane separates $i$ points of $\omega_C$ from
the other $k-i$,
and so that the hyperplane is at equal distance to the closest points,
thus cutting $C$ into $k$ cuboids. Transform the cuboids by pushing the
separating walls (defined by the hyperplanes) in the direction of the
first axis until the volumes of the new cuboids are all equal. Each of
the transformed cuboids contains the image of exactly one $\xi\in
\omega
_C$ by the transformation; call this image $\mathrm{aux}(\xi,v,0)$ and the
transformed cuboid $C_{v,0}^{\omega,\xi}$. We call this stage of the
construction the \textit{initial stage}, and $C_{v,0}^{\omega,\xi}$ the
\textit{initial cell} of $\xi$.

Suppose that for every nonnegative integer $n<N$, a local allocation
between $C=v+[0,2^n)^d$ and $\omega\cap C$, if nonempty, has been
defined for any $v\in\R^d$. For a $\xi\in\omega\cap(v+[0,2^n)^d)$ let
$C_{v,n}^{\omega,\xi}$ be the cell assigned to $\xi$ by this local
allocation, and suppose that an ``auxiliary point'' $\operatorname{aux}(\xi
,v,n)\in
C_{v,n}^{\omega, \xi}$ is also given for $\xi$. Let $\Aux
_{v,n}(C):=\{
\operatorname{aux}(\xi,v,n)\dvtx   \xi\in C\cap\omega\}$. Define now a local
allocation between $C=v+[0,2^N)^d$
and $\omega\cap C$ in $d$ steps as follows.

For every $\xi\in C\cap\omega$, let $v'\in v+2^{N-1}\{0,1\}^d$ be such
that $\xi\in\omega\cap(v+[0,2^{N-1})^d)$, and consider the local
allocation between $v'+[0,2^{N-1})^d$ and $\omega\cap
(v'+[0,2^{N-1})^d)$. Then\vspace*{1pt} $\operatorname{aux}(\xi,v',N-1)$ and
$C_{v',N-1}^{\omega
,\xi}$
are the auxiliary point and the cell, respectively, for $\xi$ by this
local allocation. We will define $\operatorname{aux}(\xi,v,N)$ and
$C_{v,N}^{\omega
,\xi}$ in $d$ steps. Let $\hat C_{\xi,0}:=C_{v',N-1}^{\omega,\xi}$ and
$\auxi_{\xi,0}:=\operatorname{aux}(\xi,v',N-1)$.\vspace*{1pt}

Let $\D_i$ be the following collection of cuboids. Each element of $\D
_i$ will be a translation of the cuboid $[0,2^{N-1})^{d-i}\times
[0,2^N)^i$, and the elements of $\D_i$ are such that their disjoint
union is $C$. Formally, $\D_i= \{([0,2^{N-1})^{d-i}\times
[0,2^N)^i)+v+(k_1,\ldots, k_{d-i},0,\ldots,0)2^{N-1}\dvtx   (k_1,\ldots,
k_{d-i},0,\ldots,0)\in\{0,1\}^d \}$.
For\break $i=1,\ldots,d$, do the following. For each $K\in\D_i$, consider the
hyperplane orthogonal to the \mbox{$(d-i+1)$th} axis that splits $K$ into two
congruent parts. Let $K_1$ and $K_2$ be these two congruent parts.
Consider cuboids $K_1'$ and $K_2'$ that have the following properties:
\begin{longlist}[(iii)]
\item[(i)] they also partition $K$;
\item[(ii)] $\overline{K_1'}\cap\overline{K_2'}$ is parallel to $\overline
{K_1}\cap\overline{K_2}$ (where $\overline X$ denotes the closure of a
set $X\subset\R^d$);
\item[(iii)] ${\lambda(K_1')\over\lambda(K_2')}={\llvert  K_1\cap\omega\rrvert  \over
\llvert  K_2\cap\omega\rrvert  }$. If $\llvert  K_1\cap\omega\rrvert  =0$ ($\llvert  K_2\cap\omega\rrvert  =0$), then\vspace*{1pt}
$K_1'$ ($K_2$) is a degenerate cuboid.
\end{longlist}
One can obtain $K_1'$ from $K_1$ by an affine transformation. For each
$\xi\in K_1$ (which is equivalent to $\auxi_{\xi, i-1}\in K_1$), let
the image of $\hat C_{\xi,i-1}$ under this translation be called
$\hat C_{\xi,i}$, and let the image of $\auxi_{\xi, i-1}$ be $\auxi
_{\xi
, i}$. Proceed similarly for $K_2'$: it is the image of $K_2$ by an
affine transformation, and for $\xi\in K_2$ we define $\auxi_{\xi, i}$
and $\hat C_{\xi,i}$ as the images of $\auxi_{\xi, i-1}$ and
$\hat C_{\xi,i-1}$ under this transformation, respectively. At the end
of the cycle, $i=d$. Define $C_{v,N}^{\omega,\xi}:=\hat C_{\xi,d}$ and
$\operatorname{aux}(\xi, v, N):=\auxi_{\xi,d}$.

Let\vspace*{1pt} us fix an arbitrary $v \in\R^d$. Up to this point we only defined
the auxiliary points and cells ($\operatorname{aux}(\xi,v,N)$, $C_{v,N}^{\omega
,\xi
}$) for any $\xi\in\omega\cap(v+[0,2^{N})^d)$. Now we will extend
the definition to all configuration points $\xi\in\omega$. Define
$\operatorname{aux}(\xi,v,N)$ as the point $\operatorname{aux}(\xi,v',N)$, where $v'$
is the
unique element of $v+2^N\Z^d$ so that $\xi\in v'+[0,2^N)^d$.
Similarly, for any $v\in\R^d$ define $C_{v,N}^{\omega,\xi}$ as the
cell $C_{v',N}^{\omega,\xi}$, where $v'$ is the unique element of
$v+2^N\Z^d$ so that $\xi\in v'+[0,2^N)^d$. We remark that auxiliary
points and corresponding cells can also be defined for noncenters, as
their images by the respective affine transformations of the cells as
above. This is the point in the construction where the statement of
(\ref{e-1}) becomes evident, as $f_{v,n}^{\omega, \xi}$ will be defined
as the indicator function of $C_{v,n}^{\omega, \xi}$.

For each fixed $N$, we have defined the algorithm recursively, but we
will prefer to think about it as an algorithm running through $N$ \textit{stages}, where in the $i$th stage $C_{v,i}^{\omega,\xi}$ and $\operatorname{aux}(\xi
,v,i)$ are constructed for each $\xi\in\omega$. Each stage consists of
$d$ \textit{steps}, one for each axis.

Next we present the AKT weak allocation scheme in an infinite setup.
(This will not yet be a weak allocation rule because equivariance
fails.) Let $\omega$ now be distributed according to the Poisson point
process. For each $v \in\R^d$ and $n \geq1$ we define $\mathrm
{AKT}_{\omega,v,n}$ to be the weak allocation, whose restriction to $C$
($C \in v+2^n\Z^d+[0,2^n)^d$ with $C \cap\omega\neq\varnothing$) is
the local allocation between $C$ and $C \cap\omega$ given by the
previous method. We call the algorithm that produces $\mathrm
{AKT}_{\omega,v,n}$ AKT($v$) run up to stage $n$, or simply AKT($v$).
Note that for any $n<n'$ and $i\leq n$ the transformations taking place
in the $i$th step of the algorithm are the same for AKT($v$) run up to
stage $n$ and for AKT($v$) run up to stage $n'$. Hence the latter can
be thought of as a continuation of the former.



For a simpler discussion, condition now on $0\in\omega$, and take
$\xi=0$.

\begin{lemma} \label{egyseg}
For $\P'$-almost all $\omega' \in\Omega'$, the sequence of cells $\{
C_{v,n}^{\omega',0}\}_{n \geq1}$ of $0$ resulting from $\mathrm
{AKT}_{\omega',v,n}$
satisfies
\[
\lim_{n \to\infty} \lambda\bigl(C_{v,n}^{\omega',0}
\bigr)=\lim_{n \to
\infty} \int_{\mathbb R^d}
f_{v,n}^{\omega',0}(x) \,dx = 1,
\]
where convergence is uniform in $v$.
\end{lemma}


\begin{lemma}\label{fontos}
There exist $ c,b >0$ and an increasing family of events $(E_R)_{R>0}$
such that:
\begin{longlist}[(iii)]
\item[(i)] $1-\P'[E_R] < c\exp(-bR^d)$;
\item[(ii)] for every $\omega'\in E_R$ and $v
\in\mathbb{R}^d$, the diameter of the cell $C_{v,n}^{\omega',0}\cup
\{
0\}$ of $0$ in
the AKT($v$) run up to stage $n$ is at most $cR$;
\item[(iii)] there are constants $c_n(R)$ such that the series $\sum_{n
=1}^\infty c_n(R)$ is summable, and
for every $\omega'\in E_R$ and $v \in\R^d$ the bound
\[
\bigl\llVert  f_{v,n}^{\omega
',0}-f_{v,n+1}^{\omega',0}\bigr\rrVert  _1=\lambda\bigl(C_{v,n}^{\omega',0}\Delta
C_{v,n+1}^{\omega',0}\bigr) < c_n(R)
\]
holds for every $n$.
\end{longlist}
\end{lemma}

The proof follows from the usual analysis of the AKT algorithm, with a
slight modification needed because of uniformity in $v$; for details,
see Section~\ref{technical}.

\section{The allocation rule}\label{rule}


In this section we will provide the necessary details for the
construction that was sketched in Section~\ref{s.cons} and prove that
it is indeed well defined. Simultaneously we will execute the analysis
on the tail behavior of the diameter of certain cells in order to
verify Theorem~\ref{main}.
Recall that $f_{v,n}^{\omega,\xi}$ was defined as the indicator
function of the cell of $\xi$ resulting from
AKT($v$) run up to stage $n$, where $\xi\in\omega\in\Omega$.

The next lemma is straightforward from the definitions.

\begin{lemma}\label{konnyu}
For $\P$-almost every $\omega\in\Omega$, for almost every $v\in\R^d$
and every $n\geq0$,
the map $
f_{v,n}^{\omega}\dvtx  \omega\to L^1 (\R^d)$, $f_{v,n}^{\omega} (\xi
)=f_{v,n}^{\omega,\xi}
$ is a weak allocation.
\end{lemma}

Of course, the map $\omega\mapsto f_{v,n}^{\omega}$ in the lemma is
far from being equivariant; hence it does not define a weak allocation rule.

\begin{lemma}\label{Lemma2} For $\P'$-almost every $\omega' \in
\Omega
'$ and fixed $n\geq1$, the map $v \mapsto f_{v,n}^{\omega',0}$ is
$L^1$-continuous in $v \in[0,2^n)^d$, except for possibly $v$ in the
union of
countably many hyperplanes.
\end{lemma}

\begin{pf} 
The function $f_{v,n}^{\omega',0}$ changes continuously in $v$, except
for the case when a hyperplane $H$ of the form $H=v+\{ (x_1,\ldots,
x_{j-1}, k, x_{j+1},\ldots, x_d)\dvtx  x_i\in\R\}$, $j\in\{1,\ldots,
d\}$, $k\in\{1,\ldots, 2^n\}$ contains some point of $\omega'$. To see
this, let $v$ and $v'$ be two points such that each point of $\omega'$
is on the same side of the respective pairs of hyperplanes above.
Recall from Section~\ref{s.akt} that the constructed cells do not vary
from the viewpoint of the reference points in this case. More
precisely, the indicator function $f^{\omega',0}_{v,n}$ is the
translate of $f^{\omega',0}_{v',n}$ by $v-v'$.

The above containment property is equivalent to the condition that $v$
is in $\bigcup_{\xi\in\omega' \cap[-2^n,2^n)^d}\bigcup_{j=1}^d\bigcup_{k=1}^{2^n} -\xi+ \{ (x_1,\ldots, x_{j-1}, k, x_{j+1},\ldots,
x_d)
\dvtx  x_i\in\R^d \}$. This is a countable union of sets of measure $0$
with probability $1$; hence the points of discontinuity form a set of
measure $0$.
\end{pf}

Define\vspace*{1pt} [as in (\ref{f_n})] the averaging function
$f_n^{\omega, \xi}=\frac{1}{2^{nd}} \int_{[0,2^n)^d}
f_{v,n}^{\omega,\xi
} \,dv$ for $\omega\in\Omega$, $\xi\in\omega$. By the previous
lemma, the
integral in the definition exists.

\begin{proposition}\label{elso}
For $\P$-almost every $\omega\in\Omega$, for every $\xi\in\omega
$, the
$L^1$ limit $f^{\omega,\xi}$ of $f_n^{\omega,\xi}$ exists as $n\to
\infty
$, it is a function
with values in $[0,1]$, it has integral~1, and satisfies $\P[\diam(\{
\xi\}\cup
\supp(f^{\omega,\xi}))>R \mid  \xi\in\omega]\leq c \exp(-bR^d)$ with
some $c,b>0$.
\end{proposition}

\begin{pf}
Since the measure $\P$ is equivariant, it is enough to prove the claim
for the Palm version: for $\P'$-almost every $\omega'\in\Omega'$ and
with $\xi=0$.
Recall that $E_R$ is a monotone increasing family of events that
exhausts a subset of $\P'$-measure 1 in~$\Omega'$, and that
the function $f^{\omega',0}$, if exists, does not depend
on $R$. Fix $R$, let $E_R$ be as in Lemma~\ref{fontos} and assume that
$\omega'\in E_R$.

For an arbitrary $m, n\in\Z^+$, $m\leq n$, $u\in2^{m}\{0,1,\ldots,2^{n-m}-1\}^d$, let
\[
g_n^{m,u, \omega'}=g_n^{m,u}:=\frac{1}{2^{md}}
\int_{u+[0,2^{m})^d} f^{\omega',0}_{v,n} \,dv.
\]
In particular, $f^{\omega',0}_n=g_n^{n,u}$. Then
%
\begin{equation}
f^{\omega',0}_n=\frac{1}{2^{(n-m)d}} \sum
_{u\in2^{m}\{0,1,\ldots,2^{n-m}-1 \}^d} g_n^{m, u}.\label{e0}
\end{equation}
If ${u,u'\in2^{m}\cdot\{0,1,\ldots,2^{n-m}-1 \}^d}$, then
$g_{m}^{m,u}=g_{m}^{m,u'}$ because the sequence of dyadic partitions
used in the construction for such a $u$ and $u'$ is
the same up to stage $m$, also showing (\ref{e-1}).
We have
\begin{eqnarray*}
\bigl\llVert g_n^{m,u}-g_{m}^{m, u}\bigr
\rrVert _1 &=& \biggl\llVert 2^{-m
d}\int_{u+[0,2^{m})^d}
f^{\omega',0}_{v,n} \,dv-2^{-m d}\int_{u+[0,2^{m})^d}
f^{\omega',0}_{v,m} \,dv\biggr\rrVert _1
\\
&\leq& 2^{-m d}\int_{u+[0,2^{m})^d} \bigl\llVert
f^{\omega',0}_{v,n}-f^{\omega
',0}_{v,m} \bigr\rrVert
_1 \,dv,
\end{eqnarray*}
and similarly for $u'$.

By Lemma~\ref{fontos}(iii), the right-hand side of this inequality
is bounded from above by $\sum_{i=m}^\infty c_i(R)$.

Hence by the triangle inequality we obtain
\begin{eqnarray}\label{e1}
\qquad \bigl\llVert g_{n}^{m,u}-g_{n}^{m,u'}
\bigr\rrVert _1 &\leq& \bigl\llVert g_n^{m,
u}-g_{m}^{m, u}
\bigr\rrVert _1
+\bigl\llVert g_{m}^{m, u}-g_{m}^{m, u'}
\bigr\rrVert _1 +\bigl\llVert g_n^{m, u'}-g_{m}^{m,
u'}\bigr\rrVert _1
\nonumber\\[-8pt]\\[-8pt]\nonumber
&\leq& 2\sum_{i=m}^\infty c_{i}(R).
\nonumber
\end{eqnarray}
Using (\ref{e0}) this implies
%
\begin{eqnarray}\label{e3}
\bigl\llVert f^{\omega',0}_n-g_n^{m,0}\bigr
\rrVert _1&\leq&\sum_{u\in2^{m}\{0,1,\ldots,2^{n-m}-1 \}^d}\frac{1}{2^{(n-m)d}}
\bigl\llVert g_n^{m, u}-g_n^{m,0}\bigr
\rrVert _1
\nonumber\\[-8pt]\\[-8pt]\nonumber
&\leq& 2 \sum_{i=m}^\infty c_i (R).
\end{eqnarray}
On the other hand, for any $m,n\in\Z^+$, $m \leq n$,
\begin{eqnarray*}
\bigl\llVert f^{\omega',0}_n-f^{\omega',0}_m\bigr
\rrVert _1 &\leq&\bigl\llVert f^{\omega',0}_n-g_n^{m,0}
\bigr\rrVert _1 +\bigl\llVert g_n^{m,0}-f^{\omega',0}_m
\bigr\rrVert _1
\\
&=&\bigl\llVert f^{\omega',0}_n-g_n^{m,0}
\bigr\rrVert _1 +\bigl\llVert g_n^{m,0}-g_m^{m,0}
\bigr\rrVert _1.
\end{eqnarray*}
The first term on the right is $\leq2\sum_{i=m}^\infty c_i (R)$ by
(\ref{e3}), and the second term is $\leq\sum_{i=m}^n c_i (R)$ by
Lemma~\ref{fontos}(iii). We conclude that $(f^{\omega',0}_n)$ is a
Cauchy sequence,
and so there is a limit $f^{\omega',0}$ in $L^1$. The fact that
$f^{\omega',0}$ takes values in
$[0,1]$ follows directly from the same fact about $f^{\omega',0}_n$.
By Lemma
\ref{egyseg} it is easy to see that $\int f^{\omega',0}_n\to1$, and by
Lemma~\ref{fontos} the support of
each $f^{\omega',0}_n$ is within radius $cR
$ around 0; hence the dominated convergence theorem implies that $\int
f^{\omega',0}=1$. The bound on the tail probability of the support is
the consequence of Lemma~\ref{fontos}(i). The above hold for every
$R$ and $\P'[\cup E_R]=1$; hence the proposition follows.
\end{pf}

The next proposition implies the variant of Theorem~\ref{main} for
fractional allocation rules instead of allocation rules; see
Definition~\ref{def1}.

\begin{proposition} \label{fractional}
The map $\eta\dvtx \omega\mapsto f^{\omega}$ with $f^{\omega}\dvtx \omega
\to
L^1(\R^d), \xi\mapsto f^{\omega,\xi}$
is a fractional
allocation rule. It satisfies
\[
\P\bigl[\diam\bigl(\{0\}\cup\supp f^{\omega,0}\bigr)>R \mid 0\in \omega\bigr]
\leq c\exp \bigl(-bR^d\bigr)
\]
for some $c$ and $b >0$.
\end{proposition}

\begin{pf}
By Proposition~\ref{elso} we have that $\eta\dvtx \omega\mapsto f^{\omega}$
with $f^{\omega}\dvtx \omega\to L^1(\R^d)$, $\xi\mapsto f^{\omega',\xi}$
satisfies $(1')$ and $(3)$ in the definition of a
fractional allocation rule and similarly for the claim about the
support of
$f^{\omega,0}$. Measurability and equivariance are clear from the
construction. So it only remains to prove $(2')$.

By Lemma~\ref{konnyu}, $\xi\mapsto f_{v,n}^{\omega,\xi}$ defined on
$\omega$ is a weak allocation; in particular, (2) from Definition~\ref
{def1} holds. Hence for almost every
$x\in\R^d$,
\begin{eqnarray*}
\sum_{\xi\in\omega} f_n^{\omega,\xi}(x)&=&\sum_{\xi\in\omega}
2^{-nd}\biggl(\int_{[0,2^n)^d} f_{v,n}^{\omega,\xi} \,dv\biggr) (x)
\\
&=&
2^{-nd}\int_{[0,2^n)^d} \sum_{\xi\in\omega}f_{v,n}^{\omega,\xi}
(x) \,dv \leq
2^{-nd}\int_{[0,2^n)^d} 1 \,dv=1,
\end{eqnarray*}
showing that $\{f_n^{\omega,\xi}\dvtx   \xi\in\omega\}$ satisfies
$(2')$.

Finally, this implies
$\sum_{\xi\in\omega}
\lim_n f_n^{\omega,\xi}(x)=\lim_n \sum_{\xi\in\omega}
f_n^{\omega,\xi}(x) \leq1$ (applying the dominated convergence
theorem for every compact subset of $\R^d$).
\end{pf}
%




The following theorem is a special case of the Campbell--Mecke formula.
We will use it in the proof of Lemma~\ref{veges}.

\begin{theorem}[{\cite{SKM}}] \label{cmformula}
For any integrable $f\dvtx  \R^d \times\Omega\to\R^+$,
\[
\E\biggl[\sum_{x \in\omega} f(x, \omega)\biggr] = \int
_{\R^d} \E\bigl[f(x, \omega) \mid x \in\omega\bigr] \,\dd x,
\]
where $\E$ is expectation with respect to the Poisson point process of
unit intensity.
\end{theorem}

\begin{lemma}\label{veges}
$\P$-almost surely for almost every $x\in\R^d$, $x$ is contained in
$\supp(f^{\omega, \xi})$ for only finitely many $\xi$'s.
\end{lemma}

\begin{pf}
Let $z\in\Z^d$, and denote by $Y_z$ the random variable that is the
number of centers $\xi$ of $\omega$ such that the intersection of
$\supp(f^{\omega, \xi})$ and $z+[0,1)^d$ is nonempty. Then
\[
\E Y_z \leq\E\biggl[\sum_{\xi\in\omega} A(\xi,
\omega)\biggr],
\]
where $A(x, \omega)=1 \textrm{ if }x \in\omega$ and $\diam
(\{
x\} \cup\supp f^{\omega, x}) > \llvert  x-z\rrvert   - \sqrt{d} $, and $A(x, \omega
)=0$ otherwise.
Using Proposition~\ref{fractional} and Theorem~\ref{cmformula} we have
\begin{eqnarray*}
\E Y_z &\leq&\int_{\R^d} \P\bigl[\diam\bigl(\{x\}
\cup\supp f^{\omega,x}\bigr) > \llvert x-z\rrvert -\sqrt{d} \mid x \in\omega
\bigr]\,\dd x
\\
&\leq&\int_{\R^d} c\exp\bigl(-b\bigl(\llvert x\rrvert -
\sqrt{d}\bigr)^d\bigr) \,\dd x < \infty.
\end{eqnarray*}
Hence $\P[Y_z = \infty]=0$, and also $\P[\bigcup_{z \in\Z^d} \{
Y_z=\infty
\}]=0$, which implies the statement.
\end{pf}



A direct consequence is the following:

\begin{lemma}\label{darabok}
For $\P$-almost every $\omega$ one can partition $\R^d$ to countably many
measurable sets of the form $S=\bigcap_{i=1}^k \supp
(f^{\omega, \xi_i} )$ (with some $\xi_1,\ldots, \xi_k\in\omega$ and
$k \in\N$).
\end{lemma}

The previous lemma will enable us to define an allocation rule from our
fractional allocation rule, in such a way that the cell allocated to a
center $\xi\in\omega$ is contained in $\supp(f^{\omega,\xi})$.
Namely, for each set $S$ as in Lemma~\ref{darabok}, we will partition
$S$ into
measurable pieces $S_1,\ldots, S_k$ such that $\lambda(S_i)=\int_S
f^{\omega,\xi_i} (x)\,dx$. We will do it in a way such that $S$ and the
$f^{\omega,\xi_i}$ ($i=1,\ldots,k$) determine the
pieces $S_i$ in some previously fixed (deterministic) way, and such that
the pieces change continuously with $S$ (in terms of Hausdorff distance
between sets, say). The central issue is to obtain a partition and an
association to the centers that is translation equivariant. A method to
do so was suggested to us by Yuval Peres, replacing the original, less
elegant proof for the following lemma:

\begin{lemma}\label{final}
Let $\eta\dvtx \omega\mapsto f^{\omega}$ be a fractional
allocation rule that satisfies Lemma~\ref{darabok}. Then there is an
allocation rule $\psi\dvtx \omega\mapsto\psi_\omega$ such that for every
$\xi\in\omega$ we have $\supp( \psi_\omega(\xi)) \subset\supp(
f^{\omega,\xi})
$.
\end{lemma}

\begin{pf}
For each set of the form $S=\bigcap_{i=1}^k \supp(f^{\omega,\xi_i}
)$ as in Lemma~\ref{darabok}, let $\alpha=\{\xi_1, \dots, \xi_k\}$ and
$c_i:=\int_S (f^{\omega, \xi_i })
(x)\,dx$, $i=1,\ldots, k$. So $\sum_{i=1}^k c_i$ is the Lebesgue
measure of
$S$. Apply a version of the site-optimal Gale and Shapley algorithm
(see \cite{HHP} for a more detailed description) within $S$ to
partition it into $\mathcal S=(S_1, \dots, S_k)$ (up to a remainder set
of measure zero) as follows with $S_i$ corresponding to $\xi_i$.
First we put all the points of $S$ that are equidistant to any pair of
centers from $\alpha$ into the set $W$. 
We note that $W$ has measure zero. Now we will stage-wise define a
series of disjoint subsets of $S \setminus W$ and auxiliary sets
corresponding to the centers in $\alpha$. Let $A_1(\xi_1), \dots,
A_1(\xi_k)$ be the intersections of the Voronoi cells of $\alpha$ with
$S \setminus W$, and set $R_0(\xi_i)=\varnothing$ for all $i \in[k]$.
Suppose that we have already constructed the sets in $\{A_l(\xi
_i),R_{l-1}(\xi_i)\dvtx  i \in[k], l\leq n\}$. The disjoint sets
$A_{n+1}(\xi_1), \dots, A_{n+1}(\xi_k)$ will be obtained as follows. We
define for all $n$ and $\xi_i \in\alpha$ the rejection radius
$r_n(\xi
_i)=\inf\{r \geq0\dvtx  \lambda(A_n(\xi_i) \cap B(\xi_i,r)) \geq c_i \}$,
and the rejection sets $R_n(\xi_i)= A_n(\xi_i) \setminus B(\xi_i,
r_n(\xi_i))$. Using this notation we define $A_{n+1}(\xi_i)= (A_n(\xi
_i) \cap B(\xi_i,r_n(\xi_i))) \cup N(\xi_i)$, where $x \in N(\xi_i)$,
if and only if $x \notin\bigcup_{j=1}^n R_j(\xi_i)$, $x \in\bigcup_{j=1}^k
R_n(\xi_j)$, and if for any other $l \in[k]$, it holds that $x \notin
\bigcup_{j=1}^n R_j(\xi_l)$, then $\llvert  x - \xi_i\rrvert   < \llvert  x-\xi_l\rrvert  $.\vspace*{1pt}

With the aid of these sets we are able to define $\mathcal S$, let
$S_i=\bigcup_{n \geq1}\bigcap_{k \geq n}(A_k(\xi_i) \cap B(\xi_i,r_k(\xi
_i)))$. It is clear that the map $\xi_i \mapsto S_i$ is equivariant
under translations, that is, for any $y \in\R^d$ we have $\xi_i+y
\mapsto S_i+y$ when we run the algorithm for $S+y$ and $\alpha+y$. It
is also true that the $S_i$'s are pairwise disjoint, because the sets
$A_n(\xi_1), \dots, A_n(\xi_k)$ were also pairwise disjoint by
definition for all $n \geq1$. It remains to show that they provide a
partition of $S\setminus W$ up to a set of measure zero.
First note that for each $i$ we have $\lambda(S_i) =\liminf_{n \to
\infty} \lambda(A_n(\xi_i) \cap B(\xi_i,r_n(\xi_i)))=\limsup_{n
\to
\infty} \lambda(A_n(\xi_i) \cap B(\xi_i,r_n(\xi_i)))\leq c_i$. If there
exists an $x \notin\bigcup_{i \in[k]} S_i$, then $x \in\bigcap_{i \in
[k]}\bigcup_{n \geq1} R_n(\xi_i)$, which\vspace*{1pt} implies that $\bigcup_{n \geq1}
R_n(\xi_i)$ is nonempty for each $i \in[k]$, so therefore $\lambda
(S_i)=c_i$. It follows that $\lambda(S\setminus(\bigcup_{i \in[k]}S_i))=0$.

We also provide the intuitive picture for better understanding. Let
each $\xi_i$ start
growing a ball around itself at linear speed, simultaneously. At time $t$
let each $\xi_i$ capture all points of $S$ that its growing ball
meets, as
long as no other ball has captured the point at some earlier time.
However, when the total measure of points in $S$ that $\xi_i$ has captured
by time $t$ is equal to $c_i$, stop growing its ball; that is,
let $P_i(t')$ be equal to $P_i(t)$ as long as $t'>t$
[where $P_i(t)$
is the set of points in $S$ that have been captured by $\xi_i$ by time
$t$]. It is
clear that by some time $t$ the Lebesgue measure of $P_i(t)$
is equal to $c_i$ for each $i$. Define this $P_i (t)$ to be the part
of the cell of $\xi$ within $S$. One can also show that the
correspondence between $\alpha$ and $\mathcal S$ is stable; see \cite{HHP} for more details.

Doing this subdivision for each of the countably many
$S$, we get the cell of $\xi_i$ assigned to it by our allocation.
Measurability and invariance follow from our method and the assumptions.
\end{pf}

We finish by noting that Proposition~\ref{fractional} and Lemma~\ref
{final} imply Theorem~\ref{main}.

If we wanted the cells assigned by our allocation rule to be connected and
to contain the corresponding center, this could be done by growing
``tendrils'' that connect the pieces of each cell $C$ of the original
allocation, and the center for $C$. By taking care to preserve
measurability and equivariance, we believe that this can be done, but
we omit the details
here.

\begin{remark}
We have only worked out the allocation rule for the
translation-equivariant case. However, one can make it
isometry-equivariant. Besides parameter $v$ that determined a translate
of $2^n\Z^d$ in the definition of the function $f_{v,n}^{\omega,\xi}$,
we need to introduce a parameter $\theta\in\{x\in\R^d, \llvert  x\rrvert  =1\}$ to
determine a rotation of $\Z^d$.
When we integrate through $v$ in the definition of $f_n^{\omega,\xi}$,
we then have to integrate with respect to $\theta$ as well; otherwise
every part of the proof extends to this modified setup automatically.
\end{remark}

\begin{remark}
An intuitive interpretation of the AKT algorithm comes from thinking of
the centers as gas particles. Then
the procedure systematically equates pressures between neighboring
cubes in the dyadic subdivision.
For the sake of analysis it was easier to allow cells around particles to
expand in the directions of the axes, but a more
``canonical'' version
would be obtained without this artifact. That is, put one particle of
gas in each center of the point process, and start growing cells around
them (small balls at the beginning), whose pressures would be
proportional to the volume of the cell and at normal direction to the
surface. When two cells meet, the pressure differences between them
would tend to equate, and in the limit all the cells would have the
same volumes. The tail behavior of their diameters should be as good as
that of our fractional allocation rule (and perhaps better if we take
constant factors into account).
On an even more speculative note, we mention that the above procedure
looks like a modification of the stable allocation rule: do not fix
centers, and let the growing cells ``push'' each other while occupying
yet unoccupied territories.
\end{remark}

\begin{question} Can one make the above heuristics precise in order to
obtain a
canonical allocation rule of optimal tail?
\end{question}

\section{Proofs of Lemmas \texorpdfstring{\protect\ref{egyseg}}{2.1} and \texorpdfstring{\protect\ref{fontos}}{2.2}}\label{technical}

The cuboid $\mathrm{AKT}_{\omega,v,n}(0)$ is a result of $nd$ affine
transformations (in $n$ stages), not taking into account in how many
steps the initial cell of $0$ is constructed.
Hence we can bound the
diameter of $C_{v,n}^{\omega,0}\cup\{0\}$ by first
bounding $\llvert  \operatorname{aux}(0,v,n)\rrvert  $ and then bounding the sum of the lengths of
the edges of $C_{v,n}^{\omega,0}$ (which bounds the diameter of
$C_{v,n}^{\omega,0}$). 
In each stage there is at most $1$ step along each of the $d$ axes. Also,
the sizes of the steps along different axes are independent as random
variables. Therefore, if we wish to obtain an upper bound on the total movement
of a point $x$ during the shifts, steps along different axes can be
treated separately.

The next lemma is standard, and we prove it only for completeness.

\begin{lemma} \label{cher}
Let $X$ be a random variable with Poisson distribution of mean~$\lambda$. If $0 \leq\rho\leq2$, then
\[
\P\bigl[\llvert X-\lambda\rrvert >\lambda\rho\bigr] < 2\exp\biggl(-
\frac{\lambda\rho^2}{4}\biggr).
\]
\end{lemma}

\begin{pf} Note that the moment generating function of $X$ is
\[
M(t)=\E\bigl[\exp(tX)\bigr]=\exp\bigl(\lambda\bigl(e^t-1\bigr)\bigr),
\qquad t \in\mathbb{R}.
\]
For one side
\begin{eqnarray*}
\P[X-\lambda>\lambda\rho]&=&\P\bigl[\exp(tX)>\exp\bigl(\lambda t(\rho+1)\bigr)
\bigr]
\\
&<& \frac{\E[\exp(tX)]}{\exp(\lambda t(\rho+1))}= \frac{\exp(\lambda(e^t-1))}{\exp(\lambda t(\rho+1))},
\end{eqnarray*}
where $t > 0$, and we use Markov's inequality. Now, for $0<t<1$ we have
$e^t < 1 +t+t^2$, so
\[
\P[X-\lambda>\lambda\rho] < \frac{\E[\exp(tX)]}{\exp(\lambda
t(\rho+1))}= \exp\bigl(\lambda
\bigl(t^2+t\bigr)-\lambda t(\rho+1)\bigr).
\]
The last expression is minimized by $t=\rho/2$, so
\[
\P[X-\lambda>\lambda\rho] <\exp\biggl(-\frac{\rho^2\lambda}{4}\biggr).
\]
For the other bound,
\begin{eqnarray*}
\P[\lambda-X>\lambda\rho]&=&\P\bigl[\exp(-tX)>\exp\bigl(-\lambda t(1-\rho )\bigr)
\bigr]
\\
&<&\frac{\E[\exp(-Xt)]}{\exp(-\lambda t(1-\rho))}= \frac{\exp(\lambda(e^{-t}-1))}{\exp(-\lambda t(1-\rho))},
\end{eqnarray*}
where $t > 0$, and we use again Markov's inequality. Now, for $0<t$ we
have $e^{-t} < 1 -t+t^2/2$, so
\[
\P[\lambda-X>\lambda\rho] <\frac{\exp(\lambda(e^{-t}-1))}{\exp
(-\lambda
t(1-\rho))} = \exp\bigl(\lambda
\bigl(t^2/2-t\bigr)+\lambda t(1-\rho)\bigr).
\]
The last expression is minimized by $t=\rho$, so
\[
\P[\lambda-X>\lambda\rho] <\exp\biggl(-\frac{\rho^2\lambda}{2}\biggr).
\]\upqed
\end{pf}

For a measurable subset $B \subset\mathbb{R}^d$, let $N(B)$ denote
the number
of centers of the Poisson point process in $B$.
Let $l(C_{v,n}^{\omega,\xi}) \in\mathbb{R}^d$ denote the vector,
whose $i$th
coordinate is the length of an edge parallel to the $i$th axis of the
cuboid $C_{v,n}^{\omega,\xi}$. 

The next lemma summarizes all the needed consequences of the
concentration of the number of centers in a fixed set. Namely, the
discrepancy of this number determines the distribution of how much a
center is moved (through its auxiliary points) and a cuboid deformed
during the $\operatorname{AKT}(v)$ procedure, and these two give bounds on
the distance of the center from the resulting cell and the diameter of
the cell, respectively.

\begin{lemma} \label{shift}
There exist $ c,b >0$ and an increasing family of events $(E_R)_{R>0}$
such that:
\begin{longlist}[(iii)]
\item[(i)] $1-\P' [E_R] < c\exp(-bR^d)$;
\item[(ii)] there exist $c'_i=c'_i(R)$ ($i\in\{0,1, \ldots\}$) such that for
every $\omega'\in E_R$ and every $v\in\R^d$, one has $\llvert  \operatorname{aux}
(0,v,i)-\operatorname{aux}(0,v,i-1)\rrvert   \leq c_i'$ and $\sum_{i =0}^\infty c'_i <cR$;\vspace*{1pt}
\item[(iii)] there exist $e_i=e_i(R)$ ($i\in\{0,1,\ldots\}$) such that for
every $\omega'\in E_R$ and every $v\in\R^d$, one has
$\llvert  l(C_{v,i}^{\omega',0})-l(C_{v,i-1}^{\omega',0})\rrvert  _\infty< e_i$ and
such that $\sum_{i =0}^\infty e_i <cR$.
\end{longlist}
\end{lemma}

\begin{pf} Our analysis will loosely follow the argument of Talagrand
and Yukich \cite{TY}, although we are able to use less sophisticated
methods because the sizes of the induced displacements occurring in
further stages of the AKT algorithm decay much more rapidly in $d \geq
3$ than in $d=2$. On the other hand, we have to achieve uniform bounds
with respect to $v \in\R^d$.

It is enough to prove the lemma only for $R>R_0$ for some suitably
chosen \mbox{$R_0>0$}. Then the lemma will follow for every $R>
0$, perhaps with different constants.

Fix $R >R_0$; $R_0$ will be determined later.

In any given stage the shift of the auxiliary point of $0$ in the
direction of the $i$th axis only depends
on the number of auxiliary points in the currently considered cuboids.
Moreover, the length of the shift only depends on the
$i$th coordinates of
the auxiliary points in the cube. Hence the shifts in different
directions are independent. Therefore it will be enough to bound the
shifts along the first axis. 

We set\vspace*{1pt} $r_0(R)=\lceil\log_2 R \rceil-d-1$ for each $R>0$. That is,
$\frac{R}{2^{d+1}} \leq2^{r_0} < \frac{R}{2^d}$.

We will define the event $A_R=A_{R,1}$ in terms of bounds on the number
of centers in certain cuboids.
For each $n \geq r_0$ consider cuboids that satisfy the following three
conditions: the cuboid is the translate of $[0,2^{n-1}-2^{-n-1}) \times
[0,2^{n}-2^{-n}) \times\cdots\times[0,2^{n}-2^{-n})$, has
a corner in $2^{-n}\mathbb{Z}^d$, and either contains~$0$, or one of
its translates by $\pm(2^{n-1}-2^{-n-1},0, \ldots,0)$ does. Let
$\mathcal G^n_1$ denote the set of these objects. Similarly, for each
$n \geq r_0$ consider cuboids that satisfy the following three
conditions: the cuboid is the translate of $[0,2^{n-1}+2^{-n-1}) \times
[0,2^{n}+2^{-n}) \times\cdots\times[0,2^{n}+2^{-n})$, has
a corner in $2^{-n}\mathbb{Z}^d$ and either it contains $0$ or one of
its translates by $\pm(2^{n-1}+2^{-n-1},0, \dots,0)$ does. Let
$\mathcal G^n_2$ denote their set.
Let us set $\rho_n=\rho_n(R)=2^{-\sklfrac{5n}{4}-2d} R^{5/4}$. For
each $n \geq r_0=r_0(R)$ and $Q \in\mathcal G^n_1 \cup\mathcal
G^n_2$, define the event\vspace*{1pt} $B_Q$, that $\llvert  N(Q\setminus\{0\} )- \lambda
(Q)\rrvert   < \lambda(Q) \rho_n$. Let $A_R=\bigcap_{n=r_0}^\infty\bigcap_{Q\in
\mathcal G^n_1 \cup\mathcal G^n_2}B_Q$. Note that\vspace*{2pt} $r_0(R)$ and $\rho
_n(R)$ for any fixed $n$ are increasing in $R$, and therefore
$(A_R)_{R>0}$ is an increasing family of events.

It is straightforward from the defining formula, that $\rho_{r_0} \geq
\rho_n$ for any $n \geq r_0$. Also,
%
\begin{equation}
\label{rhobound1} \qquad\rho_{r_0}=2^{-r_0 (5/4)-2d}R^{5/4} \leq
\biggl(\frac
{2^{d+1}}{R} \biggr)^{5/4}2^{-2d}R^{5/4}=2^{(5/4)-(3d/4)}
\leq\frac{1}{2}.
\end{equation}
Furthermore,
%
\begin{equation}
\label{rhobound2} \sum_{n = r_0}^\infty
\rho_n=\rho_{r_0}\sum_{n = r_0}^\infty
\frac
{\rho
_n}{\rho_{r_0}} \leq\frac{1}{2} \sum_{i = 0}^\infty
\bigl(2^{-5/4}\bigr)^i =\frac{1}{2-2^{-1/4}}<1.
\end{equation}
First, we establish an upper bound on $(1-\P'[A_R])$. The distribution of
$N(Q\setminus\{0\})$ according to the Palm version of the Poisson
point process is Poisson with mean $\lambda(Q)$. Thus we can use Lemma
\ref{cher}
[noting that $\rho_n \leq\frac{1}{2}$ for all $n \geq r_0$ and that
$\lambda(Q)< 2^{(n-1)d}$ for $Q\in\mathcal G^n_1 \cup\mathcal
G^n_2$] and a simple union bound to get
\begin{eqnarray*}
1-\P'[A_R] &\leq&\sum_{n=r_0}^{\infty}
\sum_{Q \in\mathcal G^n_1
\cup
\mathcal G^n_2} \bigl(1-\P'[B_Q]
\bigr)
\\
&<& \sum_{n=r_0}^{\infty} \bigl(\bigl\llvert
\mathcal G^n_1\bigr\rrvert + \bigl\llvert \mathcal
G^n_2\bigr\rrvert \bigr)2 \exp\biggl(-\frac{\rho_n^2 2^{(n-1)d} }{4}
\biggr).
\end{eqnarray*}
Using $\llvert  \mathcal G^n_1\rrvert  =3(2^{2n}-1)^d$ and $\llvert  \mathcal
G^n_1\rrvert  =3(2^{2n}+1)^d$, we conclude that
\begin{eqnarray}\label{sum1}
1-\P'[A_R] &<& \sum_{n=r_0}^{\infty}2^{(2n+1)d+4}
\exp\bigl(- 2^{(n-1)d} 2^{-((5n)/2)-4d-2} R^{5/2}\bigr)
\nonumber
\nonumber\\[-8pt]\\[-8pt]\nonumber
&=&\sum_{n=r_0}^{\infty}2^{2nd+d+4} \exp
\bigl(- 2^{n(d-5/2) -5d-2} R^{5/2} \bigr).
\end{eqnarray}
Denote the $i$th term in the sum by $a_i$. Observe that
%
\begin{equation}
\label{term} a_{r_0} < R^{2d} 2^{-2d^2+d+4} \exp
\bigl(-R^d 2^{-d^2-\sklfrac{7}{2}d+\sfrac{1}{2}}\bigr),
\end{equation}
(by $R/2^{d+1}\leq2^{r_0}\leq R/2^d$), and that
\begin{eqnarray*}
\frac{a_{n+1}}{a_n}&=& 2^{2d} \exp\bigl(2^{n(d-5/2) -5d-2}
R^{\sfrac{5}{2}} \bigl(1-2^{(d-\sfrac{5}{2})}\bigr)\bigr)
\\
&\leq& 2^{2d}\exp\bigl(-2^{(n-r_0)(d-\sfrac{5}{2}) -d^2-\sklfrac{7}{2}d+\sfrac
{1}{2}} R^{d}
\bigl(2^{(d-\sfrac{5}{2})}-1\bigr)\bigr)
\\
&\leq& 2^{2d}\exp\bigl(-R^d2^{-d^2-\sklfrac{7}{2}d+\sfrac{1}{2}}
\bigl(2^{(d-\sfrac
{5}{2})}-1\bigr)\bigr).
\end{eqnarray*}

Hence there exists a constant $c'>0$ such that for $R_0$ chosen large
enough, for every $R>R_0$,
\[
1-\P'[A_R] < \exp\bigl(- c'
R^d\bigr).
\]

Now let us assume that $n\geq r_0$ throughout the following
computation. Let $W$ be an arbitrary translate of $[0,2^n)^d$
containing $0$,
with $U$ its left half and $V$ its right. Then, when conditioned on
$A_R$ (and on
$0$ being a center), the following is true:
%
\begin{equation}
\label{deviation} \biggl\llvert \frac{N(U)-N(V)}{N(U)+N(V)}\biggr\rrvert < 4
\rho_n.
\end{equation}

We can show this by making an easy observation: there are $U',V' \in
\mathcal G_1^n$ and $U'',V'' \in\mathcal G_2^n$, so that $U' \subset U
\subset U''$ and $V' \subset V \subset V''$.
Note that
%
\begin{eqnarray}\label{hanyados}
&& \frac{N(U)-N(V)}{N(U)+N(V)}
\nonumber\\[-8pt]\\[-8pt]\nonumber
&&\qquad \leq\frac{N(U \setminus\{0\})-N(V
\setminus
\{0\})}{N(U\setminus\{0\})+N(V\setminus\{0\})} +\frac
{1}{N(U\setminus
\{0\})+N(V\setminus\{0\})}.
\end{eqnarray}
On $A_R$ we have
\begin{eqnarray} \label{eq1}
&& \frac{N(U\setminus\{0\})-N(V\setminus\{0\})}{N(U\setminus
\{0\})+N(V\setminus\{0\})}\nonumber
\\
&&\qquad \leq\frac{N(U''\setminus
\{0\})-N(V'\setminus\{0\})}{N(U'\setminus
\{0\})+N(V'\setminus\{0\})}
\nonumber\\[-8pt]\\[-8pt]\nonumber
&&\qquad < \frac{1/2(2^{n}+2^{-n})^{d}(1+
\rho_n)-1/2(2^{n}-2^{-n})^{d}(1-
\rho_n)}{(2^{n}-2^{-n})^{d}(1-
\rho_n)}
\nonumber
\\
&&\qquad =\frac{1}{2} \biggl[\frac{1+\rho_n}{1-\rho_n} \biggl(\frac
{2^{n}+2^{-n}}{2^{n}-2^{-n}}
\biggr)^d-1 \biggr],\nonumber
\end{eqnarray}
where the first inequality holds by monotonicity of $N(\cdot)$ and the
second by the definition of $A_R$. To further estimate (\ref{eq1}),
use that
\[
\frac{1+\rho_n}{1-\rho_n}=1+\frac{2\rho_n}{1-\rho_n} \leq1 + 4\rho_n,
\]
(recalling $\rho_n\leq\frac{1}{2}$ for $n\geq r_0$) and that
\[
\biggl(\frac{2^{n}+2^{-n}}{2^{n}-2^{-n}} \biggr)^d= \biggl(1+\frac
{2}{2^{2n}-1}
\biggr)^d < \biggl(1+\frac{1}{2^{2n-2}} \biggr)^d < 1 +
2^d\frac{1}{2^{2n-2}}.
\]
We obtain from the two previous expressions and (\ref{eq1}) that
%
\begin{equation}
\label{eq2} \frac{N(U\setminus\{0\})-N(V\setminus\{0\})}{N(U\setminus
\{0\})+N(V\setminus\{0\})} < \frac{1}{2}\biggl(4\rho_n+
\frac{1}{2^{2n-2-d}}+\frac{1}{2^{2n-2-d}}4\rho_n\biggr),
\end{equation}
which is an upper bound for the first term on the right-hand side of
(\ref{hanyados}).
The second term on the right-hand side of (\ref{hanyados}) is (rather
roughly) bounded by
%
\begin{equation}
\label{eq3} \frac{1}{N(U\setminus\{0\})+N(V\setminus\{0\})} < \frac{1}{(2^{n}-2^{-n})^d}\frac{1}{1-\rho_n}<
\frac{1}{2^{2n}}.
\end{equation}
If $R>R_0(d)$ for an $R_0$ chosen suitably large, it follows that $\rho
_n=2^{-\vfrac{5n}{4}-2d} R^{5/4}$ is greater than $2^{-2n+2}$
for $n \geq r_0$. Therefore by adding up bounds (\ref{eq2}) and (\ref
{eq3}), we show that for $R > R_0$,
\[
\frac{N(U)-N(V)}{N(U)+N(V)} <4\rho_n.
\]
Symmetry of $A_R$ in $U$ and $V$ implies
\[
\biggl\llvert \frac{N(U)-N(V)}{N(U)+N(V)}\biggr\rrvert <4\rho_n.
\]

Having estimated the discrepancy of points in the two halves of cubes,
we are now ready to give upper bounds on the total shift during the AKT
procedure.
Let $v'$ be an arbitrary point in $\mathbb{R}^d$, and let $D_n=(\operatorname{aux}
(0,v',n))_1-(\operatorname{aux}(0,v',n-1))_1$,
the signed amount by which the auxiliary point of $0$ is translated
along the first axis in the $n$th stage of AKT($v'$).
By the observation that until the
$n$th stage, every displacement takes place inside a cube of sidelength
$2^n$, we trivially have
\[
\Biggl\llvert \sum_{n= 1}^{r_0-1}
D_n\Biggr\rrvert \leq2^{r_0-1} < \frac{R}{2^{d+1}}<
\frac{R}{2d}.
\]
%
For $n\geq r_0$, $D_n=(2^{n-1}-h_n) \frac{N(U_n)-N(V_n)}{N(U_n)+N(V_n)}$,
where $Q_n$ is the element of $v'+ 2^n \mathbb{Z}+[0,2^n)^d$ containing
$0$, with $U_n$
its left half and $V_n$ its right, and $0 \leq h_n \leq2^{n-1}$ is the
distance of $\operatorname{aux}(0,v',n-1)$ to
the hyperplane separating $U_n$ and $V_n$.
Conditioned on $A_R$ (and $0$ being a center) we then have
%
\begin{eqnarray}
\Biggl\llvert \sum_{n=r_0}^{\infty}
D_n\Biggr\rrvert &\leq&\sum_{n=r_0}^{\infty}
\biggl\llvert \bigl(2^{n-1}-h_n\bigr) \frac{N(U_n)-N(V_n)}{N(U_n)+N(V_n)}
\biggr\rrvert
\nonumber
\\
&<& \sum_{n=r_0}^{\infty} 2^{n+1}
\rho_n =\sum_{n=r_0}^{\infty}
2^{-\sfrac{n}{4}-2d+1} R^{5/4}
\nonumber\\[-8pt]\label{eq.szum0} \\[-8pt]\nonumber
&=& R^{5/4} 2^{-\sfrac{r_0}{4}-2d+1} \sum
_{n=0}^{\infty} 2^{-\sfrac{n}{4}}
\\
&\leq& R 2^{-\sklfrac{7}{4}d+\sfrac{5}{4}}\frac{1}{1-2^{-1/4}} < \frac{R}{2d}, \label{eq.szum}
\end{eqnarray}
since $\frac{R}{2^{d+1}}\leq2^{r_0}$ by the choice of $r_0$.
Thus we have that, conditioned on the event $A_R$ and $0$ being a
center, for every $v'\in\mathbb{R}^d$ the total shift of $0$ along\vspace*{1pt} the
first axis, $\llvert  (\operatorname{aux}(0,v',n))_1\rrvert  $, is at
most $\frac{R}{d}$ for every $n\geq1$
when the AKT($v'$) is run up to stage~$n$. Furthermore, we have seen in
(\ref{eq.szum0}) that on $A_R$ the length $\llvert  (\operatorname{aux}(0,v',n))_1-(\operatorname{aux}
(0,v',n-1))_1\rrvert  $ is at most $b'_n(R):=
2^{-\sfrac{n}{4}-2d+1} R^{5/4}$ whenever $n\geq r_0(R)$.

For $1 < i \leq d$ one can define the events $A_{R,i}$ similarly to
$A_{R,1}$, to establish tail
bounds of the lengths of the shifts along the $i$th axis. We define the
subevents of $A_{R,i}$ analogously to the subevents of $A_{R,1}$,
conditioning on the relative deviation (with respect to the
expectation) of the number of configuration points being between
factors $(1-\rho_n)$ and $(1+\rho_n)$ in certain cuboids, with $n$
depending on the side length of the cuboid. Thus we arrive at the same
failure probability bound as in (\ref{eq.szum}).

Now define
%
\begin{equation}
E_R:=\bigcap_{i=1}^d
A_{R,i}. \label{E_R}
\end{equation}
This event satisfies conditions (i) and (ii) of the lemma: (i)
holds by a union bound on the complement events. On the other hand,
(ii) is true by (\ref{eq.szum}) summed up for each of the $d$
directions (using the triangle inequality), with $c'_n(R)=d b'_n(R)$.
Furthermore, $(E_R)_{R>0}$ is an increasing family of events, since for
each $i$, the families $(A_{R,i})_{R>0}$ are as well.

Now we turn to the proof of (iii).

It is clear that $\llvert  l(C_{v,n}^{\omega',0})\rrvert  _\infty\leq2^{r_0-1} <
\frac
{R}{2^{d+1}}$ for $n < r_0$
because transformations of the cell only happen inside a cube of
sidelength $2^{r_0-1}$.

For $n \geq r_0$,
\[
\bigl\llvert l\bigl(C_{v,n}^{\omega',0}\bigr)\bigr\rrvert
_\infty< \bigl\llvert l\bigl(C_{v,n-1}^{\omega',0}\bigr)\bigr
\rrvert _\infty(1+ 4\rho_{n})
\]
because, conditioned on $E_R$, in each step of the $n$th stage, at most
one sidelength can be stretched to at most $(1+
4\rho_{n})$ times its current size, using (\ref{deviation}).
This implies by (\ref{rhobound2}) that
\begin{eqnarray*}
\bigl\llvert l\bigl(C_{v,n}^{\omega',0}\bigr)\bigr\rrvert
_\infty&<& \bigl\llvert l\bigl(C_{v,r_0-1}^{\omega
',0}\bigr)\bigr
\rrvert _\infty \prod_{i=r_0}^{\infty}(1+
4\rho_{i}) \leq \frac{R}{2^{d+1}} \exp\Biggl(\sum
_{i=r_0}^{\infty} 4\rho_i\Biggr)
\\
&<& \frac{R}{2^{d+1}} \exp(4).
\end{eqnarray*}
On $E_R$, each component of $l(C_{v,n+1}^{\omega',0})$ is between
multiplicative
factor $(1-4\rho_{n+1})$ and $(1+4\rho_{n+1})$ of the respective
component of $l(C_{v,n}^{\omega',0})$. This implies with the previous
computation that
\[
\bigl\llvert l\bigl(C_{v,n+1}^{\omega',0}\bigr)-l
\bigl(C_{v,n}^{\omega',0}\bigr)\bigr\rrvert _\infty< 4\rho
_{n+1}\bigl\llvert l\bigl(C_{v,n+1}^{\omega',0}\bigr)\bigr
\rrvert _\infty< 4\rho_{n+1} \frac{R}{2^{d+1}} \exp(4).
\]
The $\{\rho_n\}$ series is summable [see (\ref{rhobound2})] and $\sum_{i=r_0}^{\infty} \rho_{i} \frac{\exp(4)}{2^{d-1}} R< \frac{\exp
(4)}{2^{d-1}} R$, which proves the claim.
\end{pf}
%

%
\begin{pf*}{Proof of Lemma~\ref{egyseg}}
What we need to prove is that the measure of the cell assigned to $0$
by AKT($v$) run up to stage $n$ tends to 1 with $n$. This is again a simple
consequence of the fact that the number of Poisson points in a large
cube is
concentrated around the volume of the cube and that for any $v$ and $n$
the weak allocation $\mathrm{AKT}_{\omega',v,n}$ is the composition of
local allocations between the classes of the dyadic partition and the
centers lying in them. The only extra technicality
comes from the fact that we want to prove convergence of the cell volumes
uniformly in $v$, but this can be checked in the same way as in the
proof of
Lemma~\ref{shift}, so we only sketch it here. For $R>0$, let the events
$E_R$ be the same as in that lemma. In particular, for any $\omega'
\in
E_R$ and $v \in\R^d$ we have that $(2^n-2^{-n})^d (1-\rho_n(R))<
N(Q_n\setminus\{0\}) < (2^n+2^{-n})^d (1+\rho_n(R))$, where\vspace*{1.5pt} $Q_n$ is
the element of $v+2^n\Z^d+[0,2^n)^d$ which contains $0$. Since
$C_{v,n}^{\omega',0}$ is a cell of a local allocation between $Q_n$ and
$\omega' \cap Q_n$, $\lambda(C_{v,n}^{\omega',0}) = \frac
{2^{nd}}{N(Q_n)} \to1$ as $n \to\infty$,\vspace*{1pt} uniformly in $v$. The events
$E_R$ exhaust a set of measure $1$ in $\Omega'$, so the claim of the
lemma follows.
\end{pf*}

\begin{pf*}{Proof of Lemma~\ref{fontos}} We take the same events $E_R$
as in
Lemma~\ref{shift}, and therefore (i) is satisfied. Furthermore,
$(E_R)_{R>0}$ is an increasing family.

To show (ii), use the following upper bound together with (ii)
and (iii) of Lemma~\ref{shift}:
\[
\diam\bigl(C_{v,n}^{\omega',0}\cup\{0\}\bigr) \leq\bigl\llvert \operatorname{aux}(0,v,n)\bigr\rrvert _1 + d\bigl\llvert l\bigl(C_{v,n}^{\omega',0}
\bigr)\bigr\rrvert _\infty.
\]
To verify that (iii) holds, let us fix $n \geq r_0(R)$. Consider
the transformations
of the cell of $0$ that occur during the steps of the $n$th stage.
There are $d$
steps, and they all can be treated similarly. Therefore we only
consider the step
along the first axis.
Let $A$ denote the cell before this step, and $B$ thereafter. We
introduce an
auxiliary cell $C$: if $w \in\mathbb{R}^d$ is the shift vector of the
auxiliary point of $0$ in
this step, then let $C=A+w$. First, on $E_R$ we obtain an upper bound
on $\lambda(A \Delta C)$.
\[
\lambda(A \Delta C) \leq2 \llvert w\rrvert \bigl\llvert l(A)\bigr\rrvert
_\infty^{d-1} \leq2c_n'(R)
R^{d-1},
\]
where we use (ii) and (iii) from Lemma~\ref{shift}.
Now we bound $\lambda(C \Delta B)$ using (iii) from Lemma~\ref{shift}:
\[
\lambda(C \Delta B) \leq\bigl\llvert l(A)\bigr\rrvert _\infty^{d-1}
\bigl\llvert l(A)-l(B)\bigr\rrvert _\infty\leq R^{d-1}
e_n(R).
\]
If we consider the two series whose summands are the respective
right-hand sides of the two previous highlighted formulas, then both of
them are absolutely summable by Lemma~\ref{shift}. We obtain
$c_n(R):=d(R^{d-1}(2c'_n(R)+e_n(R)))$.
\end{pf*}


\section*{Acknowledgments}
We thank Yuval Peres and Ron Peled for discussions on gravitational
allocation and Lemma~\ref{final}.
This research was done when the second author was at the Hausdorff
Center for Mathematics in Bonn.



\printaddresses
\end{document}